\documentclass[12pt]{article}
\usepackage{amsmath,amsthm,amsfonts,amssymb,MnSymbol,mathrsfs}
\usepackage[english]{babel}

\usepackage[applemac]{inputenc}
\usepackage[all]{xy}
\numberwithin{equation}{section}

  \newcommand{\const}{\rm const}

  \newcommand{\argmax}{\rm  argmax}

  \newcommand{\grad}{\rm grad}
  \newcommand{\dom}{\rm dom}

\theoremstyle{plain}
\newtheorem{theorem}{Theorem}[section]
\theoremstyle{theorem}
\newtheorem{corollary}{Corollary}[section]
\newtheorem{lemma}{Lemma}[section]

\newtheorem{definition}{Definition}[section]
\newtheorem{remark}{Remark}[section]
\newtheorem{example}{Example}[section]

\renewenvironment{proof}{{\bf{Proof.}}}{\hfill $\Box$ \\}

\begin{document}

\title{\textbf{Asymptotic and non-asymptotic estimates for multivariate Laplace integrals}}

\footnotesize\date{}

\author{Maria Rosaria Formica ${}^{1}$, Eugeny Ostrovsky ${}^2$, Leonid Sirota ${}^2$}

\maketitle

\begin{center}
${}^{1}$  Universit\`{a} degli Studi di Napoli Parthenope, via
Generale Parisi 13,\\ Palazzo Pacanowsky, 80132,
Napoli, Italy.\\

e-mail: mara.formica@uniparthenope.it \\

\vspace{4mm}

${}^2$ Bar-Ilan University, Department  of Mathematic and Statistics, 59200 \\
Ramat Gan, Israel. \\

e-mail: eugostrovsky@list.ru\\
e-mail: sirota3@bezeqint.net \\

\end{center}

\begin{abstract}
 We derive bilateral  asymptotic  as well as non-asymptotic estimates for the multivariate Laplace integrals.

\vspace{3mm}

Possible applications:  Tauberian theorems for random vectors. \par

\end{abstract}


\vspace{2mm}


\section{Definitions. Notations. Previous results. Statement of problem.}

 \vspace{4mm}

Let $R^d,  \  d = 1,2,\ldots$ be the ordinary $d-$ dimensional
numerical (Euclidean) space
$$
R^d =  \{ x, \ x = \vec{x} = \{ x(1), x(2), \ldots, x(d)  \}  \  \},
\ \  x(j) \in R,  \ \ j = 1,2,\ldots,d,
$$
and let $(X, B, \mu)$ be a non-trivial measurable space equipped
with sigma-finite Borelian measure  $\mu$, where $X$ is a measurable
subset of $R^d$ having strictly positive measure $\mu(X) \in (0,
\infty]$. \par

\noindent Introduce the following subset of the whole  space $R^d$
$$
R^d(Z) \stackrel{def}{=} \{ x = \vec{x} \ : \ \min_i x(i) \ge Z
\}, \  \ Z = \const  \ge 1.
$$

\vspace{4mm}

\noindent {\it We will  impose in the sequel the following condition
on the set} $X$. \par

\vspace{4mm}

{\bf Condition 1.1.}
\begin{equation} \label{AZ}
\exists Z_0 = {\const} \ge 1 \ : \ \forall Z \ge Z_0 \ \Rightarrow \
X \cap R^d(Z) \ne \emptyset,
\end{equation}

 say for all the values $Z\ge 1$ sufficiently large. \par

\vspace{5mm}

Denote
$$
  \lambda = \vec{\lambda} \in R^d, \ \  (x,\lambda) =  x \cdot \lambda = \sum_{i=1}^d \lambda(i) x(i), \ \  |x|   = \sqrt{(x,x)},
 $$
 so that $ \ \dim(x) =\dim(\lambda) = d. \ $ \par

 Define also
$$
R^d_+ =  \{ x, \ x = \vec{x} = \{ x(1), x(2), \ldots, x(d)  \}  \
\}, \ \  x(j)  \ge 0, \  \ j = 1,2,\ldots,d;
$$

$$
  \Lambda(\lambda) = \Lambda = \min_i \lambda(i),  \ \  \lambda \in R^d_+;
$$
correspondingly

$$
\underline{x} = \Lambda(x) = \min_i x(i), \ \ x \in X.
$$

\vspace{3mm}

Let also  $ \  \zeta = \zeta(x), \ x \in X$, be a measurable
numerical valued continuous function $ \ \zeta: X \to R. \ $  \par

 We assume, furthermore, that $ \ \mu(X) = \infty, \ $ as long as the opposite (probabilistic) case is trivial for us. \par

\vspace{5mm}

 {\bf Definition 1.1.} {The following integral}

\begin{equation} \label{def source int}
I(\lambda) = I[\zeta](\lambda) := \int_X e^{ (\lambda, x) -
\zeta(x)  } \ \mu(dx)
\end{equation}
{is named Laplace or exponential integral.}

\vspace{5mm}

{\bf In this article we provide asymptotical as well as
 non-asymptotical upper and lower estimates of the Laplace
integral $\ I[\zeta](\lambda) = I(\lambda)$, for all sufficiently
large values of the real vector parameter $ \ \lambda =
\vec{\lambda} \in R^d_+, \ d = 1,2,3, \ldots $, say
$\Lambda(\lambda) \ge 1$ and when $\Lambda \to \infty $; we obtain
direct estimations assuming, of course, its convergence for all the
sufficiently large values of the parameter $ \ |\lambda| :=
\sqrt{(\lambda,\lambda)}$}.  \par

\vspace{4mm}

{\bf Furthermore we also obtain an inverse evaluation, i.e. we
deduce the bilateral bounds for the source function} $ \ \zeta =
\zeta(x), \ \underline{x} \ge 1, \
 \Lambda(x) \to \infty \ $, {\bf through its integral transform } $ \ I[\zeta](\lambda) \ $, {\bf with an inverse approach. } \par

\vspace{4mm}

\ The case of other \lq\lq octants\rq\rq, for instance, $ \ \lambda
\in R^d_- \stackrel{def}{=} \{\vec{\lambda}\}, \ \lambda(j) < 0 \ $
and
 $ \ \Lambda_- = \max_i \lambda(i) \to  - \infty \ $, may be investigated quite analogously. \par

\vspace{3mm}

 \ The one-dimensional case $ \ d = 1 \ $ was considered in \cite{Chen,KozOs Subg,KozOsRelations}; a preliminary result may be found in
\cite{Maslov Fedoryuk}. \par

\vspace{4mm}

 \ We will generalize the main results obtained in the articles \cite{KozOs Subg,KozOsRelations,Maslov Fedoryuk}, where are described also
 some applications of these estimates, in particular, in the probability theory.  The estimates given below may be considered in turn as a generalization of the
 classical saddle-point  method (\cite{Fedoryuk}). \par

\vspace{5mm}

 \ The paper is organized as follows. In section 2 and in section 3 we deduce respectively an upper and a lower direct estimate for the Laplace integral
 $I(\lambda)$;
section 4 and section 5 contain an investigation of the inverse
problem and, respectively, an upper and a lower estimate for the
source function through the exponential integral. In section 6 we
consider the multidimensional Tauberian  theorems for exponential
integrals; in section 7 some important examples are described. The
last section contains the concluding remarks.\par

 \vspace{5mm}

 \ Denote, as usually,

 $$
 \vec{0} = \{0,0, \ldots, 0 \}, \ \ \vec{1} = \{ 1,1,\ldots,1 \}; \  \ \dim  \vec{0} = \dim  \vec{1} = d;
 $$

$$
\vec{a} \ge \vec{b} \ \ \Leftrightarrow \ \ a(i) \ge b(i), \ \
\forall i = 1,2,\ldots,d ;
$$

$$
R^d_+(1) \stackrel{def}{=} \{x = \vec{x} \ge \vec{1} \}.
$$

 \vspace{4mm}

 \  Let us mention briefly  a possible application.
  Recall that the so-called (multivariate) moment generating function (MGF) for the random vector (r.v.) $ \ \vec{\xi} \ $ is defined by the equality
$$
 \exp \left(\phi_{\xi}(\lambda) \right) \stackrel{def}{=} {\bf E} \exp \left(\vec{\xi} \cdot \vec{\lambda} \right) =
$$

\begin{equation}
 {\bf E} \exp \left[ \ \sum_{i=1}^d \xi(i) \ \lambda(i) \ \right] = \int_{\Omega}\exp \left(\vec{\xi}(\omega) \cdot \vec{\lambda} \right) \ {\bf P}(d \ \omega) =
\end{equation}

\begin{equation}
\int_{R^d} e^{(\lambda,x)} \ f_{\xi}(x) \ dx = \int_{R^d} e^{(\lambda, x) - \ln (1/f_{\xi}(x))} \ dx,
\end{equation}
where $ \ f_{\xi}(x) \ $ denotes the density of the r.v. $ \ \xi, \ $ if there exists. \par
 \ So, the MGF function $ \  \exp \left(\phi_{\xi}(\lambda) \right) \ $
is, on the other terms, the multivariate Laplace integral.
\par

 \ It will be presumed that the r.v. $ \ \xi \ $ satisfies the so-called Cramer's  condition:

 \begin{equation} \label{Cramer}
 \exists \delta = \const > 0 \ : \forall \lambda, \ |\lambda| < \delta \ \Rightarrow \phi_{\xi}(\lambda) < \infty
 \end{equation}
and that the density function there exists. \par

\vspace{5mm}

 \ Recall that the well-known Young-Fenchel or Legendre transform for the function $ \ \zeta: X \to R \ $ is defined as follows

$$
\zeta^*(\lambda) \stackrel{def}{=} \sup_{x \in X} (\lambda \cdot x - \zeta(x)), \ \lambda \in R^d.
$$

\vspace{4mm}

 \ If some function $\phi = \phi(\lambda)$ is defined and is finite in a set $V$, i.e. $\dom[\phi] = V$, convex or not, one can define formally

$$
\phi(\lambda) = + \infty, \ \lambda \notin V,
$$
hence

$$
\phi^*(x) \stackrel{def}{=} \sup_{\vec{\lambda} \in V}  ( \vec{\lambda} \cdot\vec{x} - \phi( \vec{\lambda})), \ x \in R^d_+.
$$

\vspace{3mm}

 \ This notion plays an important role in the probability theory. Namely, let $ \ \xi = \vec{\xi} \ $ be a random vector for which

\begin{equation} \label{upp MGF}
{\bf E} \exp(\lambda \cdot \xi) \le \exp(\phi(\lambda)), \ \lambda \in R^d_+.
\end{equation}

 \ Then

\begin{equation} \label{upp Tail T}
T_{\xi}(x) \le \exp \left( \ - \phi^*(x) \ \right), \ x \in R^d_+,
\end{equation}
where $ \ T_{\xi} = T_{\xi}(x) \ $ denotes the tail function for the
r.v. $ \ \xi: \ $

\begin{equation} \label{taild1}
T_{\xi}(x) \stackrel{def}{=}{\bf P}(\vec{\xi} \ge \vec{x}), \ x \in R^d_+,
\end{equation}
the so-called  generalized Chernoff's inequality, see e.g.
\cite{Chernoff1,Chernoff2,KozOsRelations}. \par

\vspace{3mm}

 \ Moreover, this assertion may be reversed under some natural conditions (smoothness, convexity  etc.) in the following sense. Suppose $ \ d = 1 \ $
(one-dimensional case) and that the last estimate \eqref{upp Tail T}
holds true. Then, under appropriate conditions (see
\cite{KozOsRelations}),
\begin{equation} \label{upp MGF}
{\bf E} \exp(\lambda \cdot \xi) \le \exp(\phi(C_1 \cdot \lambda)), \
\ \lambda \in R^1_+
\end{equation}
for some finite constant $ \ C_1. \ $ \par

\vspace{5mm}

\section{Main result. A direct approach. Upper estimate.}

\vspace{5mm}

 \ Let us introduce some preliminary notations and conditions. Put

\begin{equation}\label{def K}
K(\epsilon) = K[X,\mu, \zeta](\epsilon) := \int_X e^{-\epsilon \ \zeta(x)} \ \mu(dx),
\end{equation}
here and in the sequel $ \ \epsilon = \const \in (0,1). \ $ \par

\vspace{4mm}

\begin{lemma}
 Assume $X = R^d_+$. Let $\mu$ be the classical Lebesgue measure and let $ \zeta = \zeta(x), \ x \in  R^d_+ $, be a non-negative
strictly convex continuous differentiable function. The function $ \
K(\epsilon), \ \epsilon > 0$, defined by \eqref{def K}, satisfies
the following estimate

\begin{equation}
K(\epsilon) \le C[\zeta,d] \ \exp(-C_0 \, \epsilon)  \,
\epsilon^{-d}, \ \  C_0 = {\const}  \in R.
\end{equation}
\end{lemma}
\vspace{4mm}

\begin{proof}
There exist  positive constants $ \ C_1, C_2, \ldots,C_d \ $ and a
number $ \ C_0 \in R \ $ such that
$$
\zeta(\vec{x}) \ge C_0 + \sum_{i=1}^d C(i) \ x(i).
$$

Indeed, one can apply the well-known Fenchel-Morau theorem

$$
\zeta(x) = \sup_{y \in R^d_+} ( \ (x,y) - \zeta^*(y)),
$$
so that, for an arbitrary $ \ y_0 \in R^d_+ \ $,

$$
\zeta(x) \ge  (x,y_0) - \zeta^*(y_0).
$$

Therefore

\begin{equation*}
\begin{split}
K(\epsilon) &\le \int_{R^d_+} \exp [ - \epsilon( \ C_0 +
\sum_{i=1}^d
C(i) \ x(i)  \ ) \ ] \ \prod_{i=1}^d d x(i) \\
 &= \left[ \
\prod_{i=1}^d C(i) \ \right]^{-1} \ e^{-C_0 \, \epsilon} \,
\epsilon^{-d}.
\end{split}
\end{equation*}

\end{proof}

Furthermore, define

$$
Z(\epsilon) = Z[X,\mu, \zeta](\epsilon) := \int_X \exp (\zeta( \ x(1 - \epsilon)) - \zeta(x) \ ) \ \mu(dx).
$$

\begin{definition}
{\rm Let $D\subset X$ be a non-empty subset of the whole set $X$. We
introduce the so-called {\it regional} Young-Fenchel
 transform for the function $ \ \zeta(\cdot) \ $

$$
\zeta^*[D](\lambda) \stackrel{def}{=} \sup_{x \in D}  (\lambda \cdot
x - \zeta(x)), \ \ \lambda \in R^d,
$$
so that

$$
\zeta^*[X](\lambda) = \zeta^*(\lambda).
$$
}
\end{definition}
\vspace{4mm}

 \ We  represent now {\it three } methods for an {\it upper} estimate of $ I(\lambda)$ for sufficiently large values of the real parameter $ \ |\lambda|. \ $ \par

{\bf A.}  \ First of all note that if the measure   $  \ \mu \ $ is
bounded: $ \  \mu(X) = M \in (0, \infty); \ $ then the integral  $ \
I(\lambda) \ $ satisfies a very simple estimate

\begin{equation}
 I(\lambda)  \le M \cdot \sup_{x \in X}  \exp \left( \lambda x - \zeta(x)    \right) =
M \cdot \exp \left( \zeta^*(\lambda)  \right).
\end{equation}

 \ Let now $   \ \mu(X) = \infty $ and let $ \   \epsilon = \const \in (0,1). \ $

 \vspace{4mm}

 {\bf B.}  \ It will be presumed the finiteness of the integral $ \ K(\epsilon) = K[X,\mu,\zeta](\epsilon) \ $
 at least for some positive value  $ \   \epsilon_0 \in (0,1)$, i.e.
$$
 K(\epsilon) < \infty   \ \  \forall \epsilon \ge \epsilon_0.
 $$

 \ It is proved in particular in   \cite{KozOsRelations}  that

\begin{equation} \label{est K}
I(\lambda) \le K(\epsilon) \cdot  \exp \left\{ (1 - \epsilon) \zeta^* \left(  \frac{\lambda}{1 - \epsilon} \right)  \right\} \le
K(\epsilon) \cdot  \exp \left\{ \ \zeta^* \left(  \frac{\lambda}{1 - \epsilon} \ \right) \ \right\}.
\end{equation}

 \ Note that in  \cite{KozOsRelations} was considered the one-dimensional case $ \ d = 1; \ $ but the general one may be investigated
 quite analogously. In detail, let $ \  \epsilon \in (0,1) \ $ be some number for which $ \ K(\epsilon) \in (0,\infty). \ $ Consider
 the following probability measure, more precisely, the family of probability measures

 $$
 \nu_{\epsilon}(A) = \frac{1}{K(\epsilon)} \ \int_A e^{ - \epsilon \zeta(x)  } \ \mu(dx),
 $$
 or symbolically
$$
 \nu_{\epsilon}(dx) = \frac{1}{K(\epsilon)} \  e^{ - \epsilon \zeta(x)  } \ \mu(dx),
 $$
 so that

 $$
 \nu_{\epsilon}(X) = \int_X \nu_{\epsilon}(dx) = 1.
 $$
We have
$$
\frac{I(\lambda)}{K(\epsilon)} = \int_X e^{ (\lambda,x) - (1 - \epsilon) \zeta(x)  } \ \nu_{\epsilon}(dx) \le
$$

$$
\exp \left\{  \sup_{x \in X} [ (\lambda,x) - (1 - \epsilon)  \zeta(x)    ]     \right\} = \exp \left\{ (1 - \epsilon) \zeta^*(\lambda/(1 - \epsilon) ) \right\}.
$$
 \ So, the relation \eqref{est K} is proved. \par

\vspace{5mm}

 \ As a  slight consequence:

\begin{equation} \label{IKeps1}
I(\lambda) \le \inf_{\epsilon \in (0,1)}
\left[ K(\epsilon) \cdot  \exp \left\{ (1 - \epsilon) \zeta^* \left(  \frac{\lambda}{1 - \epsilon} \right)  \right\} \right];
\end{equation}

\begin{equation} \label{IKeps2}
I(\lambda) \le \inf_{\epsilon \in (0,1)}
\left[ K(\epsilon) \cdot  \exp \left\{ \zeta^* \left(  \frac{\lambda}{1 - \epsilon} \right)  \right\} \right].
\end{equation}

\vspace{4mm}

{\bf C.} \ An opposite method, which was introduced in a particular
case in \cite{KozOs Subg}, \cite{KozOsRelations}.
 Define the following integral

$$
Z(\epsilon) = Z[\zeta, \mu, X](\epsilon) := \int_X e^{ \zeta((1-\epsilon) x) - \zeta(x)  } \mu(dx),
$$
if, of course, it is finite  at least for some value  $ \ \epsilon
\in (0,1).\  $. \par

Let again  $ \   \epsilon = \const \in (0, 1)$. Applying the
well-known Young inequality

 $$
 (\lambda,x) \le \zeta((1 -\epsilon) x) + \zeta^*(\lambda /(1-\epsilon)),
 $$
 we have

 $$
 I(\lambda) \le e^{ \zeta^*(\lambda/(1 -\epsilon))  } \ \int_X e^{ \zeta ((1-\epsilon) x) - \zeta(x)  } \mu(dx) = Z(\epsilon) \  e^{ \zeta^*(\lambda/(1 -\epsilon))  }.
 $$
Of course

\begin{equation}\label{estimate I inf}
I(\lambda) \le \inf_{\epsilon \in (0,1)} \left[ \ Z(\epsilon)  \ e^{  \zeta^*(\lambda/ (1 -\epsilon))  } \ \right].
\end{equation}

\vspace{4mm}

{\bf D.}  \ Denote

$$
R(\epsilon) = R[X,\mu,\zeta](\epsilon) := \min(K(\epsilon), \ Z(\epsilon)), \ \epsilon \in (0,1).
$$

 \ We conclude

\begin{equation} \label{R eps}
I(\lambda) \le  R[X,\mu,\zeta](\epsilon) \ e^{  \zeta^*(\lambda/ (1 -\epsilon))  }.
\end{equation}

\vspace{4mm}

 \  Furthermore,  we will use the following elementary inequality

$$
1 + \epsilon < \frac{1}{1 - \epsilon} \le 1 + 2 \epsilon, \ 0 < \epsilon \le 1/2.
$$

 \ Let us introduce a new function $ \ \phi(\lambda) := \zeta^*(\lambda), \ $

\begin{equation} \label{pi kappa}
 \pi_{\kappa}(\lambda) \stackrel{def}{=} \frac{\kappa}{(\lambda \cdot \zeta^{*'}(\lambda))} = \frac{\kappa}{(\lambda \cdot \phi^{'}(\lambda))}, \  \ \kappa = \const \in (0, \infty),
\end{equation}

\begin{equation} \label{pi 1}
 \pi(\lambda) = \pi_1(\lambda) \stackrel{def}{=} \frac{1}{(\lambda \cdot \zeta^{*'}(\lambda))} = \frac{1}{(\lambda \cdot \phi^{'}(\lambda))} ,
\end{equation}
alike ones in the monograph \cite{Os Mono}, chapter 3; and suppose that

\begin{equation} \label{pi to 0}
\lim_{\Lambda(\lambda) \to \infty} \pi(\lambda) = 0;
\end{equation}
so that the value $ \ \Lambda_0 = \Lambda(\kappa)\ $ may be chosen
such that

$$
\forall \lambda: \ \min_i \lambda(i) \ge \Lambda(\kappa) \ \
\Rightarrow  \ \ \pi_{\kappa}(\lambda) \le 1/2.
$$
Let us impose the following condition on the function $ \
\phi(\cdot): \ $

\begin{equation} \label{condition phi}
\sup_{\min_i \lambda(i) \ge \Lambda(\kappa)} \left[ \phi \left( \lambda + \frac{2 \lambda}{(\lambda, \phi'(\lambda))} \ \right)  - \phi(\lambda)\ \right] = C(\phi,\kappa) = C(\phi) < \infty.
\end{equation}

\vspace{4mm}

\noindent Define also
$$
r(\lambda)= r[\zeta,\kappa](\lambda) \stackrel{def}{=} R[X,\mu,\zeta](\pi_{\kappa}(\lambda)).
$$

\vspace{4mm}

\noindent Choosing $ \ \epsilon = \pi(\lambda) \ $  in the domain $
\ \min_i \lambda(i) \ge \Lambda(1) \ $, we have the following
\begin{theorem}\label{upper th 2.1}
If the function $\phi(\cdot)=\zeta^*(\cdot)$ satisfies the condition
\eqref{condition phi}, then
\begin{equation} \label{upp gen}
I(\lambda) \le e^{C(\phi)} \ r(\lambda) \ e^{ \zeta^*(\lambda) }.
\end{equation}
\end{theorem}

\vspace{4mm}

\begin{example}
{\rm Assume in addition $ \  R(\epsilon) \le C_1 \
\epsilon^{-\beta}, \ \beta = \const \in (0,\infty), \ \epsilon \in
(0,1); \  $ then

$$
I(\lambda) \le C_1 \ e^{C(\phi)} \ \epsilon^{-\beta} \ e^{2 \epsilon/\pi(\lambda)} \ e^{\phi(\lambda)}, \ \lambda(i) \ge \Lambda,
$$
and, after the minimization over $ \ \epsilon \ $,

$$
I(\lambda) \le C_1 \ e^{C(\phi)} \ 2^{\beta} \ \beta^{-\beta} \ e^{\beta} \ \pi_1^{-\beta}(\lambda) \ e^{\zeta^*(\lambda)}.
$$
}
\end{example}
\vspace{5mm}

{\bf E.} Let us consider  an arbitrary  simple partition $ \ X = X_0 \cup X_1, \  X_0 \cap X_1 = \emptyset  \ $ of the whole set $ \ X \ $
onto two disjoint measurable subsets. We deduce splitting integral $ \ I(\lambda) \ $ into two ones

$$
I(\lambda) = \int_{X_0} \exp(\lambda x - \zeta(x)) \ \mu(dx) + \int_{X_1} \exp(\lambda x - \zeta(x)) \ \mu(dx) = I_0 + I_1,
$$
and  applying the foregoing estimates:

$$
I_0 \le \mu(X_0) \exp \left[ \ \sup_{x \in X_0}(\lambda x - \zeta(x)) \ \right] = \mu(X_0) \ \exp \left[ \ \zeta^*[X_0](\lambda) \ \right],
$$

\begin{equation*}
\begin{split}
I_1 & \le R[X_1,\mu,\zeta] (\epsilon) \ \exp \left[ \ \sup_{x \in
X_1} (\lambda x/(1 - \epsilon) - \zeta(x)) \  \right] \\
& =R[X_1,\mu,\zeta] (\epsilon) \ \exp \left[
\zeta^*[X_1](\lambda/(1 - \epsilon)  \right].
\end{split}
\end{equation*}
Denote
$$
 W[X,\mu,\zeta, \epsilon](\lambda) = :\mu(X_0) \ \exp \left[ \ \zeta^*[X_0](\lambda) \ \right] +  R[X_1,\mu,\zeta] (\epsilon) \ \exp \left[  \zeta^*[X_1](\lambda/(1 - \epsilon)  \right],
$$

\begin{equation}
W_0[X,\mu,\zeta](\lambda) = \inf_{\epsilon \in (0,1)}  W[X,\mu,\zeta, \epsilon](\lambda).
\end{equation}

\vspace{5mm}

\noindent We obtained actually  the following compound estimate.

\vspace{2mm}

\begin{lemma}
Suppose
\begin{equation}
\exists c \in (0,\infty) \ : \  \forall \lambda, \ |\lambda| \ge c \
\Rightarrow \  W_0[X,\mu,\zeta](\lambda) < \infty.
\end{equation}
Then, $ \ \forall \lambda: |\lambda| \ge c \ $ and $ \ \forall
\epsilon \in (0,1)$,

\begin{equation}
I(\lambda) \le    W[X,\mu,\zeta, \epsilon](\lambda).
\end{equation}
\end{lemma}
As a slight consequence:
\begin{equation}
I(\lambda) \le    W_0[X,\mu,\zeta](\lambda), \  \ |\lambda| \ge c.
\end{equation}

\vspace{5mm}

\begin{remark}
{\rm Introduce the following condition on the function $ \
\zeta(\cdot)$:

\begin{equation}
\exists C_1 \in [1,\infty) \ : \  W_0[X,\mu,\zeta](\lambda) \le \
\exp \left\{ \ \zeta^*(C_1 \ \lambda) \ \right\}, \ \ |\lambda| \ge
c.
\end{equation}

 \ This condition is satisfied if, for example, the function $ \ \zeta = \zeta(x), \ x \in X \ $ is regular varying:
\begin{equation} \label{slowly1}
\zeta(\lambda) = |\lambda|^m \ L(|\lambda|), \ \  |\lambda| \ge 1,
\end{equation}
where $ \ m = \const > 0, \ |\cdot| \ $ is the ordinary Euclidean
norm (or an arbitrary other non-degenerate vector one) and $ \  L =
L(r), \ r \ge 1$, is some positive continuous slowly varying
function as $ \ r \to \infty$, and we suppose

$$
\forall A \in B \ \Rightarrow \mu(A) = \int_A |x|^{\alpha} \ M(|x|)
\ dx, \ \  \alpha = {\const} > -d,
$$
where, as before,  $ \  M = M(r), \ r \ge 1$, is some positive
continuous slowly varying function as $ \ r \to \infty$. Briefly: $
\ \mu(dx) = |x|^{\alpha} \ M(|x|) \ dx.  \  $ We have $ \
K(\epsilon) \le \overline{K}(\epsilon), \ \epsilon \in (0,1);\ $

$$
\overline{K}(\epsilon) := \int_{R^d} \exp \left( \  - \epsilon |x|^m L(|x|) \  \right) \ |x|^{\alpha} \ M(|x|) \ dx.
$$
One can apply the spherical coordinates:

 $$
 \overline{K}(\epsilon) =  \frac{\pi^{d/2}}{\Gamma(d/2 + 1)} \, K_0(\epsilon),
 $$
where
$$
K_0(\epsilon)=\int_0^{\infty} \exp \left(  -\epsilon \ r^m \ L(r)
\right) \ r^{\alpha + d - 1} \ M(r) \ dr.
$$
We obtain, after the substitution $ \ r^m \epsilon  = y,  \  dr =
m^{-1} y^{1/m - 1} \ \epsilon^{-1/m} \,dy$,
\begin{equation*}
\begin{split}
Z_m(\epsilon) & \stackrel{def}{=} \ m \ \epsilon^{ (\alpha + d)/m
}\ K_0(\epsilon) \\
&=\int_0^{\infty} y^{ (\alpha + d)/m - 1  } \exp
\left( \ -y \ L(y^{1/m} \epsilon^{-1/m}) \  \right) \ M
\left(y^{1/m} \ \epsilon^{-1/m} \right) \ dy
\end{split}
\end{equation*}
and, as $ \ \epsilon \to 0+ \ $,
\begin{equation*}
\begin{split}
Z_m(\epsilon) &\sim M(\epsilon^{-1/m}) \int_0^{\infty} e^{-y \
L(\epsilon^{-1/m}) } \ y^{ (\alpha + d)/m - 1 } \ dy \\
&= M(\epsilon^{-1/m}) \ \Gamma((\alpha + d)/m) \ L^{ -(\alpha +
 d)/m}(\epsilon^{-1/m}),
\end{split}
\end{equation*}
where $\Gamma$ is the classical Gamma function.

\vspace{4mm}

To summarize: as $ \ \epsilon \to 0+ \ $

\begin{equation} \label{overline K}
   m \ \overline{K}(\epsilon) \sim
\frac{\pi^{d/2}}{\Gamma(d/2 + 1)} \ \epsilon^{-(\alpha + d)/m} \ \Gamma((\alpha + d)/m) \ \frac{M \left(\epsilon^{-1/m} \right)}{L^{(\alpha+d)/m} \left(\epsilon^{-1/m} \right)}.
\end{equation}

\vspace{4mm}

 Thus, in this case, the values $K = K(\epsilon)$
and $R = R(\epsilon)$, \ $\epsilon \in (0,1)$, are finite with
concrete estimate following from \eqref{overline K}:
\begin{equation} \label{overline Z}
   m \ Z(\epsilon) \le C[\zeta,m,d] \
\frac{\pi^{d/2}}{\Gamma(d/2 + 1)} \ \epsilon^{-(\alpha + d)/m} \ \Gamma((\alpha + d)/m) \ \frac{M \left(\epsilon^{-1/m} \right)}{L^{(\alpha+d)/m} \left(\epsilon^{-1/m} \right)}.
\end{equation}
}
\end{remark}

\vspace{5mm}

If the condition of Remark 2.1 is satisfied, then

\begin{equation}
I(\lambda) \le  \exp \left( \  \zeta^*(C_2 \ \lambda) \   \right),  \ L(\lambda) \ge 1.
\end{equation}

\vspace{6mm}

\begin{theorem}
Let $ \  X = R^d_+   \ $ and $ \ \mu \ $ be the ordinary Lebesgue
measure. Suppose that the random vector $ \ \vec{\xi} \ $, with
non-negative entries $  \ \{\xi(i) \}, \ i = 1,2,\ldots,d \  $,
satisfies the Cramer's condition:
$$
\exists  \  \lambda_0 = \vec{\lambda_0} = \{ \lambda_0(i) \}, \ i =
1,2,\ldots, \ \lambda_0(i) > 0  \ : \ {\bf E} \exp \left( \lambda_0
\cdot \xi  \right) < \infty.
$$
Then
$$
 \exists \epsilon_0 > 0 \ : \ \forall \epsilon > \epsilon_0 \ \Rightarrow \ K(\epsilon) = K[X,\mu,\zeta](\epsilon) < \infty.
$$
\end{theorem}

\vspace{4mm}

\begin{proof}
Denote for brevity  $ \   G(x) =  G[\xi](x), \ $  so that

$$
T_{\xi}(\vec{x}) = e^{-G(\vec{x})}, \ x \ge 0.
$$

\noindent It is sufficient to consider only the two-dimensional
case: assume

$$
B  = B(\lambda,\mu):= \int_0^{\infty} \int_0^{\infty} e^{\lambda \ x + \mu \ y - G(x,y) } \ dx dy  < \infty
$$
for some positive values $ \ \lambda, \mu. \ $ We have

$$
B =  \sum_{n=0}^{\infty} \sum_{m=0}^{\infty} \int_n^{n+1} \int_m^{m+1}   e^{\lambda \ x + \mu \ y - G(x,y) } \ dx dy  \ge
$$

$$
 \sum_{n=0}^{\infty} \sum_{m=0}^{\infty} \int_n^{n+1} \int_m^{m+1} e^{ \lambda n + \mu m  - G(n+1, m+1) } \ dx dy =
 \sum_{n=0}^{\infty} \sum_{m=0}^{\infty} e^{ \lambda n + \mu m  - G(n+1, m+1) },
$$
therefore

$$
\sum_{n=0}^{\infty} \sum_{m=0}^{\infty} e^{ \lambda n + \mu m  - G(n+1, m+1) } < B(\lambda,\mu) < \infty,
$$
so

$$
e^{ \lambda n + \mu m  - G(n+1, m+1) } \le B(\lambda,\mu) < \infty, \ G(n+1,m+1) \le B \ e^{-\lambda n - \mu m},
$$
and finally

$$
(\lambda \ \mu)^{-1} K(\epsilon) < \int_0^{\infty} \int_0^{\infty} e^{-\epsilon G(x,y)} dx \ dy \le \sum_{n=0}^{\infty} \sum_{m=0}^{\infty}
\exp \left\{ -\epsilon[\lambda (n-1) + \mu(m-1)] \right\}  < \infty,
$$
if

$$
\epsilon > \epsilon_0 := \max \left(\lambda^{-1}, \mu^{-1} \right),
$$
\end{proof}

\vspace{4mm}

\section{Main result. A direct approach.  Lower estimate.}

 \vspace{4mm}

 We introduce additional notations.

$$
S(\lambda,x) = (\lambda,x) - \zeta(x), \ x_0 = x_0(\lambda) \in
{\argmax}_{x \in X} S(\lambda,x),
$$
where, by definition,
$$
{\argmax}_{x \in X} S(\lambda,x) = \{ x \in X  \ : \ S(\lambda, \ x)
= \zeta^*(\lambda) \  \}.
$$

Obviously, the value $ \ x_0 = x_0(\lambda)  \ $ may be non-unique.
\par

Furthermore, we introduce the variables
$$
X_0 = X_0(\epsilon) = X_0(\epsilon, \lambda) := \{x \in X \ : \
S(\lambda,x) \ge \zeta^*(\lambda(1 - \epsilon)) \}, \ \epsilon =
{\const} \in (0,1),
$$

$$
U(\epsilon) = U[\zeta](\epsilon,\lambda) := \int_{X_0(\epsilon)} \mu(dx) = \mu(X_0(\epsilon,\lambda)).
$$

\vspace{4mm}

\begin{theorem}\label{lower estimate I-U}
Let $\epsilon \in (0,1)$ be such that $U(\epsilon)>0$. Then, for
sufficiently large values $ \min_i \lambda(i) \ge \Lambda = {\const}
\ge 1$, we have
$$
I(\lambda) \ge U[\zeta](\epsilon,\lambda) \ \exp \left(
\zeta^*(\lambda(1 - \epsilon) )  \right), \ \epsilon \in (0,1), \ \
\min_i \lambda(i)  \ge \Lambda.
$$
Of course,
$$
I(\lambda) \ge \sup_{\epsilon \in (0,1)} [ \
U[\zeta](\epsilon,\lambda) \ \exp \left(  \zeta^*(\lambda(1 -
\epsilon) )  \right) \ ], \ \  \min_i \lambda(i)  \ge \Lambda.
$$
\end{theorem}

\vspace{4mm}

\begin{proof}

\begin{eqnarray*}
\begin{split}
I(\lambda) & = \int_X \exp \left[ \ \lambda \ x - \zeta(x) \ \right]
\mu(dx)  \ge \int_{X_0} \exp [\ \lambda \ x - \zeta(x) \ ] \mu(dx)
\\
&\ge \int_{X_0} \exp \left[ \ \zeta^*(\lambda(1 - \epsilon)) \
\right] \ \mu(dx)  =  U[\zeta](\epsilon,\lambda) \ \exp \left[ \
\zeta^*(\lambda(1 - \epsilon) ) \ \right].
\end{split}
\end{eqnarray*}
\end{proof}

As a slight consequence we get:

\begin{corollary}

$$
I(\lambda) \ge U[\zeta](\epsilon,\lambda) \ \exp \left( \  \zeta^*(\lambda) -  \epsilon \ (\lambda, \zeta^{*'}(\lambda)  )  \ \right),
$$
and, if we choose $ \ \epsilon = \pi_{\kappa}(\lambda), \ $
\begin{equation}\label{estimate Cor3.1}
I(\lambda) \ge U[\zeta](\pi_{\kappa}(\lambda),\lambda) \ \exp \left( \  \zeta^*(\lambda) - \kappa \ \right).
\end{equation}
\end{corollary}
\noindent Let us define the following function
\begin{equation} \label{estim V}
V(\lambda) = V[\zeta](\lambda) \stackrel{def}{=} \sup_{\kappa > 0}
\left\{ \ U[\zeta](\pi_{\kappa}(\lambda), \lambda) \ e^{-\kappa} \
\right\},
\end{equation}
so, by \eqref{estimate Cor3.1}, we have
\begin{equation} \label{Low estim}
I(\lambda) \ge V[\zeta](\lambda) \ e^{ \ \zeta^*(\lambda) \,},  \ \
\min_i \lambda(i) \ge \Lambda.
\end{equation}

\vspace{5mm}

For instance, it is reasonable to suppose in addition, see e.g.
Example 3.1 below, that

$$
 U[\zeta](\pi_{\kappa}(\lambda),\lambda) \ge \gamma \ \Lambda^{\alpha} \ \kappa^{\beta} \ , \ \alpha,\beta,\gamma = \const \in (0,\infty);
$$
then

$$
I(\lambda) \ge \gamma \ (\beta/e)^{\beta} \ \Lambda^{\alpha} \ e^{ \zeta^*(\lambda) }, \ \Lambda = \min_i \lambda(i) \ge e.
$$

\vspace{5mm}

Let us consider the following important example. \par

 \vspace{3mm}

\begin{example}
{\rm
 Suppose that $ \  X = R^d_+, \ d \mu = dx$ and that
 the function $ \ \zeta = \zeta(x), \ x \in X = R^d_+ \ $ is non-negative, strictly convex,
 twice continuous and differentiable as well as its conjugate  $ \ \zeta^*(\lambda) \ $ and such that its second (matrix) derivative

$$
\zeta^{''}(x) = \left\{ \frac{\partial^2 \zeta(x)}{\partial x(i) \ \partial x(j)}   \right\}, \ i,j = 1,2,\ldots,d
$$
is a strictly positive definite matrix for all sufficiently large
values $ \ \min_i x(i). \ $ \par

Denote also
$$
\zeta'(x) = {\grad} \zeta(x) = \left\{ \ \frac{\partial \zeta}{
\partial x(i)} \ \right\},
$$

$$
x_0 = \vec{x}_0[\zeta](\lambda) = x_0[\zeta](\lambda) = {\argmax}_{x
\in R^d_+} S(\lambda,x) =  {\argmax}_{x \in R^d_+} [(\lambda,x) =
(\lambda,x) - \zeta(x)],
$$

$$
\Delta = \Delta(\lambda,x) = S(\lambda,x)  - S(\lambda(1 -
\epsilon), x_0(\lambda)),
$$
so that

$$
{\grad} \zeta(x_0) = \lambda,  \ \ \lim_{\Lambda \to \infty}
x_0[\zeta](\lambda) = \infty
$$
and
\begin{equation*}
\begin{split}
X_0(\epsilon,\lambda) &= \{ x \in R^d_+ \ : \ S(\lambda,x) \ge
\zeta^*(\lambda(1 - \epsilon))\\
 &=\{\ x \in R^d_+ \ : \
S(\lambda,x) \ge S(\lambda(1 - \epsilon), x_0(\lambda)) \  \}.
\end{split}
\end{equation*}

We deduce after simple calculations, using Taylor's formula, that
the set \ $ X_0(\epsilon,\lambda) \ $  is asymptitical equivalent,
as $ \  \epsilon \to 0+ \ $, to the following one (multidimensional
ellipsoid)
$$
\tilde{X}_0 = \left\{ x: (\zeta^{''}(x_0)(x-x_0), (x - x_0)\ ) \le
\epsilon \ (\lambda, x_0)  \ \right\},
$$
in the sense that
$$
\lim_{\epsilon \to 0+} \frac{\mu(\tilde{X}_0)}{\mu(X_0)} = 1.
$$
The case when the value $\epsilon=\epsilon(\lambda)$ is dependent on
$\lambda$, but such that
$$
 \lim_{\Lambda \to \infty}  \epsilon(\lambda) = 0,
$$
can not be excluded. \par

It is no hard to calculate the "volume" of ellipsoid \ $
\tilde{X}_0:\ $

$$
\mu \left(\tilde{X}_0 \right)  = \frac{\pi^{d/2} \ \left[ 2 \epsilon \cdot(\lambda, x_0) \right]^d \ }{\Gamma(d/2 + 1)} \ \cdot \left\{ \ \det \zeta^{''}(x_0) \ \right\}^{-1/2}.
$$
Following

$$
\mu \left(X_0 \right) \ge  C_0(d) \ [ \ \epsilon \cdot (\lambda, x_0[\zeta](\lambda)  \ ]^d \cdot \left\{ \ \det \zeta^{''}(x_0) \ \right\}^{-1/2}.
$$

\vspace{4mm}

If, for instance, $ \ d = 1, \ m = {\const} > 1, \ $

$$
  \zeta(x) = \zeta_m(x) \stackrel{def}{=} x^m/m, \
$$
then

$$
  \ x_0 = \lambda^{1/(m-1)}, \ \  \lambda, x \ge 1; \ \ \mu(X_0) \ge C_m \epsilon^{1/2} \ \lambda^{1/(m-1)},
$$
and we find, after some calculations,

\begin{equation*}
\begin{split}
I_m(\lambda) & \stackrel{def}{=} \int_{0}^{\infty}\exp(\lambda x -
x^m/m) \ dx
\\
&\ge \left( \frac{2.5}{m-1} \right)^{1/2} \ \lambda^{(2- m)/(2m -
2)} \ \exp \left(\lambda^{m'}/m' \right),
\end{split}
\end{equation*}
where  $m' = m/(m-1), \ \lambda \ge \lambda_0(m)$ and
\begin{equation*}
\lambda_0(m)=\left\{
  \begin{array}{ll}
    2 & \hbox{if} \ \  1< m \leq 2 \\
     \displaystyle 2 \left( \frac{m-2}{2m - 2}  \right)^{(m-1)/m} & \hbox{if} \ \  m> 2
  \end{array}
\right.
\end{equation*}

The last estimate is in full accordance, up to a multiplicative
constant, with the exact asymptotic estimates for $I_m(\lambda)$, as
$\lambda \to \infty$,  which may be find, e.g., in the well-known
book \cite{Fedoryuk}, sections 1, 2:

\begin{equation}
I_m(\lambda) \sim   \sqrt{2\pi/(m-1)} \ \lambda^{(2- m)/(2m - 2)} \ \exp \left(\lambda^{m'}/m' \right).
\end{equation}

\vspace{3mm}

The upper estimate corresponding to the lower one obtained above,
for the integral $I_m(\lambda)$, has the form

$$
I_m(\lambda) \le m^{2/m - 1} \ e^{1/m} \ \Gamma(1/m) \ \lambda^{1/(m-1)} \ \exp \left( \lambda^{m'}/m'  \right), \ \lambda \ge \lambda_0(m).
$$
}
\end{example}

\vspace{5mm}

\begin{center}
\section{Inverse approach. Upper estimation.}
\end{center}

\vspace{5mm}

Let now the representation \eqref{def source int} be given on the
form of an {\it inequality}

\begin{equation}\label{upper J}
 \int_X e^{  (\lambda, x) -  \zeta(x)  } \ \mu(dx) \ge J(\lambda), \  \ \lambda \in R^d_+
\end{equation}
for a certain  non-negative continuous function $ J =  J(\lambda)$.
Here we derive the upper bound for the source function $ \ \zeta =
\zeta(x) \ $ for all the sufficiently large values
 $ \ \Lambda(x) = \min_i x(i), \ i = 1,2,\ldots,d$, of course, under appropriate conditions. \par

Let us impose the following condition on our datum. Namely, assume
that for some finite constant $ \ C_{12} \ $

\begin{equation}
e^{C(\zeta^*)}  \ r(\lambda) \ e^{\zeta^*(\lambda)} \le e^{ \ \zeta^*( \ C_{12}[\zeta](\lambda) \ ) \ }, \ \lambda \ge \Lambda = \const \ge 1,
\end{equation}
Suppose also that the function $ \ \zeta(\cdot) \ $ is non-negative,
continuous and convex. We have, by virtue of Theorem \ref{upper th
2.1},

$$
e^{ \ \zeta^*( \ C_{12}[\zeta](\lambda)\ ) \ } \ge J(\lambda), \ \lambda \ge \Lambda,
$$

$$
\ln J(\lambda) \le \zeta^*(C_{12} \lambda), \ C_{12} = C_{12}[\zeta] \in (0, \infty),
$$
therefore

$$
[ \ \ln J(\cdot) \ ]^*(x) \ge \zeta^{**}( \ x/C_{12} \ ).
$$

\vspace{3mm}

 Under the above conditions and by virtue of
Fenchel-Moreau Theorem, we have

\vspace{2mm}
\begin{theorem}
If the function $\zeta(\cdot)$ satisfies the condition \eqref{upper
J}, where $ J =  J(\lambda)$ is a non-negative function, then
\begin{equation} \label{inver stat}
\zeta(x) \le [ \  \ln J(\cdot) \ ]^*(C_{12} \ x), \  \ \Lambda(x) =
\min_i x(i) \ge 1.
\end{equation}
\end{theorem}

\vspace{5mm}

\section{Inverse approach. Lower estimation.}

\vspace{5mm}

Let now the representation \eqref{def source int} be given on the
form of an {\it inequality}

\begin{equation}\label{lower K}
 \int_X e^{  (\lambda, x) -  \zeta(x)  } \ \mu(dx) \le K(\lambda), \ \lambda \in R^d_+
\end{equation}
for a certain non-negative continuous function $ \ K =  K(\lambda)$.
Here we derive the lower bound for the source function $ \ \zeta =
\zeta(x)$.
\par

Let us impose the following condition on our datum. Namely, assume
that there exists $C_{13} = C_{13}[\zeta] = \const  \in (0,1)$ such
that

\begin{equation}
e^{C(\zeta^*)}  \ V[\zeta](\lambda) \ e^{\zeta^*(\lambda)} \ge e^{ \
\zeta^*( \ C_{13} \cdot \lambda \ ) \ }, \ \lambda \ge \Lambda =
\const \ge 1.
\end{equation}

Suppose, as above, that the function $ \ \zeta(\cdot) \ $ is
non-negative, continuous and convex. We have, by virtue of Theorem
\ref{lower estimate I-U} and its consequences,

$$
e^{ \ \zeta^*( \ C_{13} \lambda \ ) \ } \le K(\lambda), \  \ \lambda
\ge \Lambda,
$$

$$
\ln K(\lambda) \ge \zeta^*(C_{13} \lambda),
$$
therefore, under the above conditions and by virtue of
Fenchel-Moreau Theorem, we have

\vspace{2mm}

\begin{theorem}
If the function $\zeta(\cdot)$ satisfies the condition \eqref{lower
K}, where $K = K(\lambda)$ is a non-negative continuous function,
then
\begin{equation} \label{inver stat}
\zeta(x) \ge [ \  \ln K(\cdot) \ ]^*(C_{13} \ x), \ \Lambda(x) \ge
1.
\end{equation}
\end{theorem}

\vspace{5mm}

\section{Multivariate Tauberian theorems.}

 \ {\bf  Preface. } Tauberian theorems are named the relations between asymptotical or not-asymptotical behavior of some function (sequence) and
correspondent behavior of its certain integral transform, for
example, Laplace, Fourier or power series transform, see
\cite{Tauber,Korevaar}. They play a very important role, for
example, in the probability theory (see \cite{Bingham}), to
establish the connection between the behavior of tail of
distribution for a random variable and the asymptotic one of its
Moment Generation Function (MGF).\par

There are many results in this direction for one-dimensional case,
as well as asymptotical ones, see e.g. in
\cite{Bagdasarov,Broniatowski,Davies,Eichel,Geluk1,Geluk2,Kasahara1,Kasahara2,Kosugi,Mason,OsSirInver,Yakimiv}.
\par

 In this section we investigate multivariate Tauberian theorems describing relations between the function $ \ \zeta = \zeta(x), \ x \in X \ $ and its Laplace integral
 transform $ \ I[\zeta](\lambda),\  \lambda \in R^d, \ $ when $ \ \Lambda(x) \to \infty \ $ or correspondingly $ \ \Lambda(\lambda) \to \infty. \ $  \par

\vspace{4mm}

\begin{center}

{\sc Direct approach.}

\end{center}

\vspace{5mm}

\begin{theorem}\label{upper limit}
{\rm (\textbf{Upper limit}).}
 Under the same assumptions of Theorem \ref{upper th 2.1} for the function $\phi(\lambda) =
\zeta^*(\lambda)$, if in addition suppose that
\begin{equation} \label{limphi}
\lim_{\min_i \lambda(i) \to \infty}  \phi(\lambda) = \infty
\end{equation}
and
\begin{equation} \label{lim r lambda}
\lim_{\min_i \lambda(i) \to \infty} \frac{|\ln
r(\lambda)|}{\phi(\lambda)} = 0,
\end{equation}
then
\begin{equation} \label{underline1}
 \overline{\lim}_{\min_i \lambda(i) \to \infty}
\frac{\ln I(\lambda)}{\phi(\lambda)} \le 1.
\end{equation}
\end{theorem}

 \vspace{4mm}

\begin{proof}
 Choosing  $ \ \epsilon = \epsilon(\lambda) = \pi(\lambda) = \pi_1(\lambda) \ $ we have, for sufficiently large values
 $ \ \Lambda(\lambda) = \min_i \lambda(i)$,

 $$
 \frac{ \ln I(\lambda)}{\phi(\lambda)} \le \frac{|\ln r(\lambda)|}{\phi(\lambda)} + \frac{\phi(\lambda + \lambda \pi(\lambda))}{\phi(\lambda)}.
 $$
 The term on the left hand side tends to zero as $\Lambda \to \infty$, the limit of the quantity on the right hand side is equal to one. In detail,

$$
\frac{\phi(\lambda + \lambda \pi(\lambda))}{\phi(\lambda)} \ge \frac{\phi(\lambda )}{\phi(\lambda)} = 1.
$$
On the other hand, from the condition \eqref{condition phi} it
follows

$$
\phi(\lambda + 2 \lambda \ \pi(\lambda)) \le C(\phi) + \phi(\lambda),
$$
therefore, by virtue of condition \eqref{limphi},

$$
\frac{\phi(\lambda + \lambda \pi(\lambda))}{\phi(\lambda)} \le 1 + \frac{C_1(\phi)}{\phi(\lambda)} \to 1,
$$
as $ \ \Lambda(\lambda) \to \infty.$ \par

This completes the proof.
\end{proof}

\vspace{5mm}

 \begin{theorem}\label{lower limit}
 {\rm (\textbf{Lower limit}).}
 Under the same assumptions of Theorem \ref{lower estimate I-U} for the function $\phi(\lambda) =
\zeta^*(\lambda)$, if in addition suppose that
\begin{equation}
\lim_{\Lambda(\lambda) \to \infty} \frac{ \ \ln \{ \ V(\lambda) \} \ }{\zeta^*(\lambda)} = 0,
\end{equation}
where $V$ is defined in \eqref{estim V},
 then
\begin{equation}
\underline{\lim}_{\Lambda (\lambda) \to \infty} \frac{ \ \ln  I(\lambda) \ }{\zeta^*(\lambda)} \ge 1.
\end{equation}
\end{theorem}
\vspace{4mm}

\begin{proof}
The proof is completely alike to the one based on Theorem \ref{upper
limit} and may be omitted.
\end{proof}

\vspace{4mm}

As consequence of Theorems \ref{upper limit} and \ref{lower limit}
we have

\begin{theorem}\label{consequence upper-lower limit}
Suppose that all the conditions of Theorems \ref{upper limit} and
\ref{lower limit} are satisfied. Then there exists the following
limit and
\begin{equation} \label{Dir lim}
\lim_{ \Lambda(\lambda) \to \infty} \frac{ \ \ln  I(\lambda) \ }{\zeta^*(\lambda)} = 1.
\end{equation}
\end{theorem}

\vspace{5mm}

\begin{center}

{\sc Inverse approach.} \par

\end{center}

\ Given the representation \eqref{def source int}, in which the
function $ \ \zeta = \zeta(x) \ $ is {\it convex and continuous}, we
have

\begin{theorem}\label{inv lower limit}
{\rm (\textbf{Lower estimate})} Suppose that there exists a
continuous non-negative function $Q = Q(\lambda), \ \lambda \in
R^d$, for which
$$
\lim_{ \Lambda(\lambda) \to \infty  } \frac{\ln V(\lambda)}{Q(\lambda)} = 0,
$$
$$
\overline{\lim}_{\Lambda(\lambda) \to \infty} \frac{\ln I(\lambda)}{Q(\lambda)} \le 1,
$$
and
$$
 \lim_{\Lambda(\lambda) \to \infty} Q(\lambda) = \infty.
$$
Then
$$
\underline{\lim}_{\Lambda(x)  \to  \infty} \frac{\zeta(x)}{Q^*(x)} \ge  1.
$$

\end{theorem}

\vspace{4mm}

\begin{proof}
 Let $\delta \in (0,1)$ be an arbitrary \lq\lq small\rq\rq number. There exists a value $ \ \Lambda_0 = \Lambda_0(\delta) > 1 \ $ such
that, for all the values $\lambda$,
$$
\Lambda(\lambda) \ge \Lambda_0 \ \ \Rightarrow \ \ \ln I(\lambda)
\le (1 + \delta) \ Q(\lambda),
$$
$$
 I(\lambda) \le \exp ( (1 + \delta) \ Q(\lambda)).
$$
We apply the estimation of Corollary \ref{estimate Cor3.1}, so that

$$
V[\zeta](\lambda) \ e^{ \ \zeta^*(\lambda) \  } \le \exp ( (1 + \delta) \ Q(\lambda)),
$$
and
$$
\frac{\zeta^*(\lambda)}{Q(\lambda)} \le \frac{\ln
V(\lambda)}{Q(\lambda)} + (1 + \delta) \le (1 + 2 \delta), \ \
\Lambda(\lambda) \ge 2 \Lambda_0.
$$
Therefore
$$
\zeta^*(\lambda)  \le ( 1 + 2 \delta) Q(\lambda),
$$
and
$$
\zeta^{**}(x) \ge \frac{1}{1 + 2 \delta} \ Q^* \left(  \frac{x}{1 + 2 \delta}  \right).
$$
Applying the Fenchel-Morau Theorem, we conclude the proof.
\end{proof}

\vspace{4mm}

\ Given again the representation \eqref{def source int}, in which
the function $ \ \zeta = \zeta(\lambda) \ $ is {\it convex and
continuous}, we have

\begin{theorem}\label{inv upper limit}
{\rm (\textbf{Upper estimate})} Suppose that there exists a
continuous non-negative function $Q = Q(\lambda), \ \lambda \in
R^d$, for which
$$
\lim_{ \Lambda(\lambda) \to \infty  } \frac{\ln r(\lambda)}{Q(\lambda)} = 0,
$$
$$
\overline{\lim}_{\Lambda(\lambda) \to \infty} \frac{\ln I(\lambda)}{Q(\lambda)} \ge 1,
$$
and
$$
 \lim_{\Lambda(\lambda) \to \infty} Q(\lambda) = \infty.
$$
Then
$$
\overline{\lim}_{\Lambda(x)  \to  \infty} \frac{\zeta(x)}{Q^*(x)} \le  1.
$$

\end{theorem}

\vspace{4mm}

\begin{proof}
The proof is quite alike as the one in Theorem \ref{inv lower
limit}. Let $\delta \in (0,1)$ be an arbitrary \lq\lq small\rq\rq
number. There exists a value $ \Lambda_0 = \Lambda_0(\delta) > 1$
such that, for all the values $ \lambda$,
$$
 \Lambda(\lambda) \ge \Lambda_0 \ \Rightarrow \ \ln I(\lambda) \ge
(1 - \delta) \ Q(\lambda),
$$
$$
 I(\lambda) \ge \exp ( (1 - \delta) \ Q(\lambda)).
$$
We apply the estimation of Corollary \ref{estimate Cor3.1}, so that

$$
e^{C(\phi) \ }r[\zeta](\lambda) \ e^{ \ \zeta^*(\lambda) \  } \ge \exp ( (1 - \delta) \ Q(\lambda)),
$$
and
$$
\frac{\zeta^*(\lambda)}{Q(\lambda)} \ge \frac{C(\phi) + \ln r(\lambda)}{Q(\lambda)} + (1 - \delta) \ge (1 - 2 \delta), \ \Lambda(\lambda) \ge 2 \Lambda_0.
$$
Therefore
$$
\zeta^{**}(x) \le \frac{1}{1 - 2 \delta} \ Q^* \left(  \frac{x}{1 - 2 \delta}  \right).
$$
Applying the Fenchel-Morau Theorem, we conclude the proof.
\end{proof}

\vspace{4mm}

 \ To summarize. \par

 \vspace{5mm}

\begin{theorem}
{\rm (\textbf{Hybrid  estimate}).} Suppose that all the conditions
of Theorems \ref{inv lower limit} and \ref{inv lower limit} are
satisfied. Then the following limit there exists and
$$
\lim_{\Lambda(x)  \to  \infty} \frac{\zeta(x)}{Q^*(x)} =  1.
$$
\end{theorem}

 \vspace{5mm}

\section{An important example.}

\vspace{5mm}

 \ In this section we consider $ \ X = R^d \ $ as well as $ \ \lambda \in R^d.   \ $ \par

\vspace{3mm}

\begin{definition}
{\rm Recall that the function $ \  g = g(x): R^d \to R \  $ is said
to be {\it radial}, or equally {\it spherical invariant}, iff it
depends only on the Euclidean norm $ \ |x| \ $ of the vector $ \ x =
\vec{x}$, namely there exists $g_0: R \to R$ such that
$$
g(x) = g_0(|x|).
$$
}
\end{definition}
\vspace{4mm}

\begin{lemma}
 Suppose that the function $ \ g: R^d \to R \ $ is radial and such that its Young-Fenchel transformation $ \ g^*(y) \ $ there exists.
Then it is again a radial function, namely there is a function $g_0:
R \to R $  for which
\begin{equation} \label{radial}
g^*(y) = g_0^*(|y|) = \sup_{z \in R} (|y| \ z - g_0(z)).
\end{equation}
 As a consequence, it is an even function.

Moreover, the optimal value in the definition of the Young-Fenchel
transformation, i.e. the variable
$$
x(y) = x[g](y) := {\argmax}_{x \in R^d}( (x,y) - g(x) ),
$$
so that $ \ g^*(y) = (y,x[g](y) - g(x[g])(y)), $ is also a radial
function if, of course, there exists and is uniquely determined.
\end{lemma}

\vspace{3mm}

\begin{proof}
Let $ \ U: R^d \to R^d \ $ be an arbitrary linear {\it unitary}
operator and let $ \ U^* \ $ be its conjugate (linear) operator,
also unitary. Recall that a function $ \ f: R^d \to R \ $ is radial
iff for an arbitrary linear unitary operator $U$, it is $ f(Ux) =
f(x), \ x \in R^d$.

We have
\begin{equation*}
\begin{split}
g^*(U y)  & = \sup_{x \in R^d} ((x,Uy) - g(x))  = \sup_{x \in R^d}
((U^*x,y) - g(x)) \\
&= \sup_{x \in R^d} ((U^*x,y) - g(U^*x))  =   \sup_{z \in R^d}
((z,y) - g(z))  =  g^*(y).
\end{split}
\end{equation*}
Therefore the function $g^*(y)$ is radial.

\noindent The second proposition has an alike proof.
\end{proof}

\vspace{4mm}

\begin{remark}
{\rm The radiality of the Fourier transform of a radial function is
well-known, see e.g. \cite{SteinWeiss}, chapters 2,3.

\vspace{5mm}

Let us consider the following family of Young-Fenchel functions

$$
\zeta_{\kappa,L}(\lambda) \stackrel{def}{=} \kappa^{-1} \ |\lambda|^{\kappa} \ L^{1/\theta} \left(|\lambda|^{\kappa} \right), \ |\lambda| \ge e,
$$

$$
\zeta_{\kappa,L}(\lambda) = C \ \lambda^2, \ |\lambda| < e.
$$
where $ \ \kappa = {\const} > 1, \theta = \kappa/(\kappa - 1), \
L(r), \ r \ge e \ $  is a slowly varying at infinity, twice
continuous and differentiable function, such that
$$
\lim_{r \to \infty} \frac{L(r)}{L(r/L(r))} = 1.
$$
The Young-Fenchel transformation for these functions is calculated
in particular in the monograph \cite{Seneta}, chapter 1, sections
1,3,4: as $ \ x \to \infty \ $
$$
\zeta_{\kappa,L}^*(x) \sim \theta^{-1} \ |x|^{\theta} \ L^{1/\theta}(x).
$$
}
\end{remark}
\vspace{3mm}

One can apply our theory of Tauberian theorems.

 \vspace{4mm}

\begin{theorem}

Denote
\begin{equation} \label{Ex case}
I_{\kappa, L}(\lambda):= \int_X e^{  (\lambda, x) -  \zeta_{\kappa,L}^*(x)  } \ dx.
\end{equation}
We have
\begin{equation} \label{ex kappa1}
\lim_{\Lambda(\lambda) \to \infty} \frac{\ln I_{\kappa, L}(\lambda) }{\zeta_{\kappa,L}(\lambda)} = 1.
\end{equation}
Furthermore, the inverse conclusion holds true. Namely, if for some
Young-Orlicz function $ \ \zeta = \zeta(x) \ $
\begin{equation} \label{ex kappa2}
\lim_{\Lambda(\lambda) \to \infty} \frac{\ln I[\zeta](\lambda)(|\lambda|) }{\zeta_{\kappa,L}(|\lambda|)} = 1,
\end{equation}
then
\begin{equation} \label{ex kappa2}
\lim_{\Lambda(x) \to \infty} \frac{\zeta(|x|)}{\zeta^*_{\kappa,L}(|x|)} = 1.
\end{equation}
\end{theorem}

\vspace{5mm}

A particular case:

 $$
 \zeta(x) = \zeta_{m,r}(x)  = m^{-1} \ |x|^m \ \ln^r(|x|), \  \ |x|\ge e, \ m = {\const} > 1, \ r = {\const} \in R.
 $$
  We obtain, after some calculations, as $|y| \to \infty$,

$$
\zeta_{m,r}^*(y) \sim  \frac{1}{m'} \ (m-1)^{r/(m-1)} \ |y|^{m'} \
[\ln |y| ]^{-r/(m-1)},
$$
where $ m' = m/(m-1)$.

\vspace{5mm}

\section{  \ Concluding remarks.}

 \vspace{3mm}

\ {\bf A.} It is interesting, by our opinion, to generalize the
estimates obtained in the second section to the case of
infinite-dimensional linear spaces, as well as to generalize our
estimates for the more general integrals of the form

$$
I[\zeta](\lambda) := \int_X \exp \zeta(\lambda,x) \ \mu(dx).
$$

\vspace{3mm}

 \ {\bf B.} One can consider also the applications of the obtained results in the Probability theory, namely, in the theory of great deviation,
 asymptotical or not.

\vspace{5mm}


\vspace{6mm}


 \begin{thebibliography}{55}


\vspace{5mm}

\bibitem{Bagdasarov}
{D.R. Bagdasarov and E.I. Ostrovsky}, {\it Reversion of Chebyshev's
Inequality}. Probab. Theory Appl., \textbf {40} (4) (1996),
737--742.

\bibitem{BennShar}
{C. Bennet and R. Sharpley},  {\it  Interpolation of operators}.
Academic Press, Inc., Boston, MA, 1988.

\bibitem{Bingham}
{N. H. Bingham}, {\it Tauberian theorems and large deviations}.
arXiv:0712.3410v1 [math.PR]  20 Dec 2007.


\bibitem{Broniatowski}
{M. Broniatowski and A. Fuchs}, {\it Tauberian Theorems, Chernoff
Inequality and the Tail Behavior of Finite Convolution of
Distribution Function}. Adv. Math., \textbf {116} (1) (1995),
12--33.


\bibitem{BuldKozMetric}
{V.V. Buldygin and Yu.V.  Kozachenko},  {\it Metric Characterization
of Random Variables and Random Processes}. Translations of
Mathematics Monograph, AMS, v.188 (1998).

\bibitem{Chen}
{H. Chen}, {\it Evaluation of the Laplace integral}. Internat. J.
Math. Ed. Sci. Tech., \textbf {35} (5) (2004), 773--777.


\bibitem{Chernoff1}
 {H. Chernoff}, {\it A career in statistics.} In X. Lin, C. Genest, D.L. Banks, G. Molenberghs, D.W. Scott,
 J-L. Wang, {\it Past, Present, and Future of Statistical Science}. CRC Press. p. 35. ISBN 9781482204964 (2014).


\bibitem{Chernoff2}
 {H. Chernoff}, {\it A measure of asymptotic efficiency for tests of a hypothesis based on the sum of
observations}.
 Ann. Math. Statistics , \textbf {23,} (1952), 493--507.


\bibitem{Davies}
{P.L. Davies}, {\it Tail probabilities for positive random variables
with entire characteristic functions of very regular growth}. Z.
Angew. Math. Mech., \textbf {56,} (1976), 334--336.


\bibitem{Eichel}
{P. Eichelsbacher and L. Knichel}, {\it Fine asymptotics for models with Gamma type moments.} \\
arXiv:1710.06484v1 [math.PR] 17 Oct 2017.


\bibitem{Fedoryuk}
{M.V. Fedoryuk}, {\it The saddle-point method}. Moscow, Nauka,
(1977) (In Russian).


%
%
%
%


\bibitem{Geluk1}
{J.L. Geluk, L. de Haan and  U. Stadtm\"{u}ller}, {\it A Tauberian
theorem of exponential type}. Canad. J. Math. \textbf {38} (3)
(1986), 697-718.

\bibitem{Geluk2}
{J.L. Geluk}, {\it On the relation between the tail probability and
the moments of a random variable}. Nederl. Akad. Wetensch. Indag.
Math. \textbf{46} (4) (1984), 401--405.


\bibitem{Janson1}
{S. Janson}, {\it Further examples with moments of Gamma
type}.\\
 arXiv:1204.5637v2, 6 Feb 2013.

\bibitem{Janson2}
{S. Janson}, {\it Moments of Gamma type and the Brownian supremum
process area}. Probab. Surv., \textbf {7} (2010), 1--52.


\bibitem{Kasahara1}
{Y. Kasahara}, {\it Tauberian theorems of exponential type}. J.
Math. Kyoto Univ. \textbf {18} (2) (1978), 209--219.

\bibitem{Kasahara2}
{Y. Kasahara and N. Kosugi}, {\it Remarks on Tauberian theorem of
exponential type and Fenchel-Legendre transform}. Osaka J. Math.,
\textbf {39} (3) (2002), 613--619.

\bibitem{Kolokoltsov}
{V.N. Kolokoltsov, T.M. Lapinski}, {\it Multivariate Laplace
approximation with estimated error and application to limit
theorems}. arXiv:1502.03266v5 [math.PR] 17 Jul 2018.

\bibitem{Korevaar}
{J. Korevaar}, {\it Tauberian theory: a century of developments}.
Grundlehren der Mathematischen Wissenschaften, \textbf {329}
Springer-Verlag, Berlin, 2004.

\bibitem{Kosugi}
{N. Kosugi}, {\it Tauberian theorem of exponential type and its
application to multiple convolution}. J. Math. Kyoto Univ.,
\textbf{39} (2) (1999), 331--346.

\bibitem{KozOs Subg}
{Yu.V. Kozachenko and E.I. Ostrovsky}, {\it The Banach Spaces of
random Variables of subgaussian Type}. Theory of Probab. and Math.
Stat. (in Russian). Kiev, KSU, 32, (1985). 43--57.

\bibitem{KozOsRelations}
{Yu.V. Kozachenko, E.I. Ostrovsky and L. Sirota}, {\it Relations
between exponential tails, moments and moment generating functions
for random variables and vectors}. arXiv:1701.01901v1 [math.FA] 8
Jan 2017.

\bibitem{Maslov Fedoryuk}
{V.P. Maslov and M.V. Fedoryuk}, {\it  Logarithmic Asymptotic
behavior of the Laplace integrals}. Mathematical Notes, \textbf{30}
(5) (1981), 763--768.

\bibitem{Mason}
{D.M. Mason}, {\it An extended version of the Erd\"os-R\'enyi strong
law of large numbers}. Ann. Probab., \textbf{17} (1) (1989),
252--265.

\bibitem{Os Mono}
{E.I. Ostrovsky}, {\it Exponential estimations for Random Fields and
its applications,} (in Russian). Moscow - Obninsk, OINPE (1999).

\bibitem{OsSir Vec r v}
{E. Ostrovsky and L. Sirota}, {\it Vector rearrangement invariant
Banach spaces of random variables with exponential decreasing tails of distributions}. \\
 arXiv:1510.04182v1 [math.PR] 14 Oct 2015.

\bibitem{OsSir Discontin}
{E. Ostrovsky and L. Sirota}, {\it Non-asymptotical sharp
exponential estimates for maximum distribution of discontinuous random fields}. \\
 arXiv:1510.08945v1 [math.PR] 30 Oct 2015


\bibitem{Os Support}
{E.I.  Ostrovsky}, {\it About supports of probability measures in
separable Banach spaces}. Soviet Math., Doklady, \textbf{255} (6)
(1980), 836--838, (in Russian).


\bibitem{OsSirInver}
{E. Ostrovsky and L. Sirota}, {\it Inversion of Tchebychev-Tchernov
inequality}. arXiv:1711.06896v1 [math.PR] 18 Nov 2017.


\bibitem{Seneta}
{E. Seneta}, {\it Regularly Varying Functions}. Springer-Verlag, New
York, 1976.

\bibitem{SteinWeiss}
{E. M. Stein and G. Weiss}, {\it Introduction to Fourier analysis on
Euclidean spaces.}
 Princeton University Press, Princeton, N.J., 1971.

\bibitem{Tauber}
{A. Tauber}, {\it Ein Satz aus der Theorie der unendlichen Reihen}.
Monatsh. Math. Phys. \textbf{8} (1) (1897), 273--277.

\bibitem{Yakimiv}
{A.L. Yakimiv}, {\it Probabilistic applications of Tauberian
theorems}. Modern probability and statistics, VSP, Leiden, 2005,
ISBN: 9067644374.

\bibitem{Zhang Zhou}
{A. Zhang and Y. Zhou}, {\it A Non-asymptotic, Sharp, and
User-friendly Reverse Chernoff-Cramer Bound.} arXiv:1810.09006v1
[math.PR]  21 Oct 2018.


\vspace{5mm}

 \end{thebibliography}
\end{document}